\documentclass[10pt]{article}
\usepackage{cite}
\usepackage{mathrsfs}
\usepackage{amsfonts}
\usepackage{amsmath}
\usepackage{amsfonts,amssymb,color}
\usepackage{dsfont}
\usepackage{curves}
\usepackage{mathrsfs}
\usepackage{pifont}
\usepackage{amssymb}
\allowdisplaybreaks

\numberwithin{equation}{section}

\date{}

\def\BigRoman{\uppercase\expandafter{\romannumeral\number\count 255 }}
\def\Romannumeral{\afterassignment\BigRoman\count255=}

\begin{document}
\title{Spectral radius and perfect $k$-matchings in $t$-connected graphs
%\thanks{}
}
\author{\small Quanru Pan, Sizhong Zhou\footnote{Corresponding author. E-mail address: zsz\_cumt@163.com (S. Zhou)}\\
\small  School of Science, Jiangsu University of Science and Technology,\\
\small  Zhenjiang, Jiangsu 212100, China\\
}

\maketitle
\begin{abstract}
\noindent A $k$-matching of a graph $G$ is a function $f:E(G)\rightarrow\{0,1,2,\ldots,k\}$ with $\sum\limits_{e\in E_G(v)}f(e)\leq k$ for each vertex $v$ of $G$, where $E_G(v)$ is the set of edges incident
with $v$ in $G$. A perfect $k$-matching of a graph $G$ is a $k$-matching $f$ satisfying $\sum\limits_{e\in E_G(v)}f(e)=k$ for any vertex $v$ of $G$. A fractional perfect matching of a graph $G$ is a function
$f:E(G)\rightarrow [0,1]$ satisfying $\sum\limits_{e\in E_G(v)}f(e)=1$ for any $v\in V(G)$. We denote by $\rho(G)$ the spectral radius of $G$. In this paper, we put forward a tight spectral radius condition
for a $t$-connected graph to possess a perfect $k$-matching and a tight spectral radius condition for the existence of a perfect $k$-matching in a $t$-connected graph with a fractional perfect matching.
\\
\begin{flushleft}

(2020) Mathematics Subject Classification: 05C50, 05C70
\end{flushleft}
\end{abstract}

\section{Introduction}

In this paper, we deal only with finite simple graphs without loops or multiple edges. Let $G$ be a graph with vertex set $V(G)$ and edge set $E(G)$. The order of $G$, denoted by $|V(G)|=n$, is the number
of its vertices. The complete graph of order $n$ is denoted by $K_n$. For any vertex subset $S$ of $G$, we denote by $G[S]$ the subgraph of $G$ induced by $S$, and by $G-S$ the subgraph obtained from $G$
by removing the vertices in $S$ together with the edges to which the vertices in $S$ are incident. The degree of a vertex $v$ in $G$, denoted by $d_G(v)$, is the number of edges incident with the vertex
$v$ in $G$. A vertex $v$ with $d_G(v)=0$ is called an isolated vertex of $G$. The set of isolated vertices in $G$ is denoted by $Iso(G)$. We write $i(G)=|Iso(G)|$. Let $odd(G)$ denote the number of odd
components with at least three vertices in $G$. Let $G_1$ and $G_2$ be two disjoint graphs. We use $G_1\cup G_2$ to denote the disjoint union of $G_1$ and $G_2$. The join of $G_1$ and $G_2$, denoted by
$G_1\vee G_2$, is the graph obtained from $G_1\cup G_2$ by adding edges joining each vertex of $G_1$ to each vertex of $G_2$.

Given a graph $G$ with $V(G)=\{v_1,v_2,\ldots,v_n\}$, the adjacency matrix of $G$, denoted by $A(G)$, is a $(0,1)$-matrix $A(G)=(a_{ij})_{n\times n}$ with $a_{ij}=1$ if $v_iv_j\in E(G)$, and $a_{ij}=0$
otherwise. The adjacency spectral radius (or spectral radius, for short) of $G$, denoted by $\rho(G)$, is the largest eigenvalue of its adjacency matrix $A(G)$. For some properties on spectral radius, we
refer the reader to \cite{PB,ELW,BB,Wc,Ws1,Ws2,ZBW,ZBS,Za,ZZZL}.

Lu and Wang \cite{LW} introduced the concept of perfect $k$-matching. A $k$-matching of a graph $G$ is a function $f:E(G)\rightarrow\{0,1,2,\ldots,k\}$ with $\sum\limits_{e\in E_G(v)}f(e)\leq k$ for each
vertex $v$ of $G$, where $E_G(v)$ is the set of edges incident with $v$ in $G$. A perfect $k$-matching of a graph $G$ is a $k$-matching $f$ satisfying $\sum\limits_{e\in E_G(v)}f(e)=k$ for any vertex $v$
of $G$. In particular, a perfect 1-matching is the perfect matching (or 1-factor) in usual sense.

O \cite{Os} put forward a spectral radius condition for the existence of a perfect 1-matching in a graph. Zhao, Huang and Wang \cite{ZHW} provided a new spectral condition for a graph to possess a perfect
1-matching. Fan, Goryainov, Huang and Lin \cite{FGHL}, Zhou \cite{Zs} established some spectral radius conditions for the existence of a perfect 1-matching in a bipartite graph. Zhang \cite{Zhang} arose an
eigenvalue condition to ensure that a regular graph contains a perfect $k$-matching. Zhang and Fan \cite{ZF} showed a lower bound on the spectral radius for a graph to possess a perfect $k$-matching.

\medskip

\noindent{\textbf{Theorem 1.1}} (Zhang and Fan \cite{ZF}). Let $k$ be an odd integer with $k\geq3$, and let $G$ be a connected graph of even order $n\geq6$. Then the following statements hold.

(\romannumeral1) If $n=6$ and $\rho(G)\geq\rho(K_2\vee4K_1)$, then $G$ contains a perfect $k$-matching, unless $G=K_2\vee4K_1$.

(\romannumeral2) If $n\geq8$ and $\rho(G)\geq\rho(K_1\vee(K_{n-3}\cup2K_1))$, then $G$ contains a perfect $k$-matching, unless $G=K_1\vee(K_{n-3}\cup2K_1)$.

\medskip

One can check that $\rho(K_1\vee(K_{n-3}\cup2K_1))>\rho(K_t\vee(K_{n-2t-1}\cup(t+1)K_1))$ for $n\geq8$ and $t\geq2$. If we consider a $t$-connected graph $G$ in Theorem 1.1, then the bound on $\rho(G)$ can
be improved to $\rho(G)\geq\rho(K_t\vee(K_{n-2t-1}\cup(t+1)K_1))$.

\medskip

\noindent{\textbf{Theorem 1.2.}} Let $t$ and $k$ be two positive integer with $k\equiv1$ (mod 2), and let $G$ be a $t$-connected graph of even order $n\geq5t+3$. If
$$
\rho(G)\geq\rho(K_t\vee(K_{n-2t-1}\cup(t+1)K_1)),
$$
then $G$ has a perfect $k$-matching, unless $G=K_t\vee(K_{n-2t-1}\cup(t+1)K_1)$.

\medskip

If $t=1$ in Theorem 1.2, then we immediately get Theorem 1.1 (\romannumeral2). Hence, Theorem 1.2 is an improvement and generalization of Theorem 1.1.

A fractional perfect matching of a graph $G$ is a function $f:E(G)\rightarrow [0,1]$ satisfying $\sum\limits_{e\in E_G(v)}f(e)=1$ for any $v\in V(G)$. Pan and Liu \cite{PL}, Li, Miao and Zhang \cite{LMZ},
Chen and Guo \cite{CG} proposed some spectral radius conditions for a graph to possess a fractional perfect matching. Xue, Zhai and Shu \cite{XZS} established a lower bound for the spectral radius in a
graph with minimum degree $\delta$ to guarantee the existence of a fractional perfect matching. Zhou \cite{Zs2,Zt} provided two spectral conditions for a graph to contain a fractional perfect matching
with given properties.

Jia, Fan and Liu \cite{JFL} showed a tight spectral radius condition for a graph with a fractional perfect matching to possess a perfect 1-matching.

\medskip

\noindent{\textbf{Theorem 1.3}} (Jia, Fan and Liu \cite{JFL}). Let $G$ be a connected graph of even order $n\geq12$ with a fractional perfect matching. If
$$
\rho(G)\geq\rho(K_1\vee(K_{n-5}\cup K_3\cup K_1)),
$$
then $G$ has a perfect 1-matching, unless $G=K_1\vee(K_{n-5}\cup K_3\cup K_1)$.

\medskip

Motivated by Theorem 1.3, it is natural and interesting to propose a tight spectral radius condition to guarantee that a $t$-connected graph with a fractional perfect matching has a perfect $k$-matching.

\medskip

\noindent{\textbf{Theorem 1.4.}} Let $t$ and $k$ be two positive integers with $k\equiv1$ (mod 2), and let $G$ be a $t$-connected graph of even order $n\geq5t+7$ with a fractional perfect matching. If
$$
\rho(G)\geq\rho(K_t\vee(K_{n-2t-3}\cup K_3\cup tK_1)),
$$
then $G$ has a perfect $k$-matching, unless $G=K_t\vee(K_{n-2t-3}\cup K_3\cup tK_1)$.

\medskip

If $t=1$ and $k=1$ in Theorem 1.4, then we immediately get Theorem 1.3. Hence, Theorem 1.4 is a generalization of Theorem 1.3. Our paper is organized as follows. In Section 2, we provide some preliminary
results. In Section 3, we verify Theorem 1.2. In Section 4, we prove Theorem 1.4.

\section{Preliminary lemmas}

In this section, we provide several preliminary lemmas, which will be used to verify our main results. Lu and Wang \cite{LW} proposed a characterization for a graph with a perfect $k$-matching.

\medskip

\noindent{\textbf{Lemma 2.1}} (Lu and Wang \cite{LW}). Let $k\geq1$ be an odd integer. A graph $G$ has a perfect $k$-matching if and only if
$$
odd(G-S)+k\cdot i(G-S)\leq k|S|
$$
for every subset $S\subseteq V(G)$.

\medskip

The following lemma provides a necessary and sufficient condition to ensure that a graph contains a fractional perfect matching.

\medskip

\noindent{\textbf{Lemma 2.2}} (Scheinerman and Ullman \cite{SU}). A graph $G$ has a fractional perfect matching if and only if
$$
i(G-S)\leq|S|
$$
for any subset $S\subseteq V(G)$.

\medskip

\noindent{\textbf{Lemma 2.3}} (Li and Feng \cite{LF}). Let $H$ be a subgraph of a connected graph $G$. Then
$$
\rho(G)\geq\rho(H),
$$
where the equality follows if and only if $G=H$.

\medskip

\noindent{\textbf{Lemma 2.4}} (Miao, Li and Wei \cite{MLW}). Let $n_1,n_2,\ldots,n_t,p$ be positive integers with $n_1\geq n_2\geq\cdots\geq n_t\geq p$ and $\sum\limits_{i=1}^{t}n_i=n-s$. Then
$$
\rho(K_s\vee(K_{n_1}\cup K_{n_2}\cup\cdots\cup K_{n_t}))\leq\rho(K_s\vee(K_{n-s-p(t-1)}\cup(t-1)K_p),
$$
with equality occurring if and only if $(n_1,n_2,\ldots,n_t)=(n-s-p(t-1),p,\ldots,p)$.

\medskip

\noindent{\textbf{Lemma 2.5.}} Graph $K_t\vee(K_{n-2t-1}\cup(t+1)K_1)$ has no perfect $k$-matching, where $t,k$ and $n$ are positive integers with $k\equiv1$ (mod 2), $n\equiv0$ (mod 2) and $n\geq2t+2$.

\noindent{\it Proof.} Let $G=K_t\vee(K_{n-2t-1}\cup(t+1)K_1)$. Set $S=V(K_t)$. If $n=2t+2$, then we have
$$
odd(G-S)+k\cdot i(G-S)=k(t+2)=kt+2k=k|S|+2k>k|S|.
$$
If $n\geq2t+4$, then we get
$$
odd(G-S)+k\cdot i(G-S)=1+k(t+1)=kt+k+1=k|S|+k+1>k|S|.
$$
Based on Lemma 2.1, graph $K_t\vee(K_{n-2t-1}\cup(t+1)K_1)$ has no perfect $k$-matching. Lemma 2.5 is proved. \hfill $\Box$

\medskip

\noindent{\textbf{Lemma 2.6}} (Fan and Lin \cite{FL}). Let $\sum\limits_{i=1}^{t}n_i=n-s$ with $s\geq1$. If $n_1\geq n_2\geq\cdots\geq n_t\geq1$ and $n_2\geq3$, then
$$
\rho(K_s\vee(K_{n_1}\cup K_{n_2}\cup\cdots\cup K_{n_t}))\leq\rho(K_s\vee(K_{n-s-t-1}\cup K_3\cup(t-2)K_1)),
$$
where the equality follows if and only if $(n_1,n_2,n_3,\ldots,n_t)=(n-s-t-1,3,1,\ldots,1)$.

\medskip

\noindent{\textbf{Lemma 2.7.}} Graph $K_t\vee(K_{n-2t-3}\cup K_3\cup tK_1)$ has no perfect $k$-matching, where $t,k$ and $n$ are positive integers with $k\equiv1$ (mod 2), $n\equiv0$ (mod 2) and $n\geq2t+6$.

\noindent{\it Proof.} Let $G=K_t\vee(K_{n-2t-3}\cup K_3\cup tK_1)$. Write $S=V(K_t)$. Then we possess
$$
odd(G-S)+k\cdot i(G-S)=2+kt=k|S|+2>k|S|.
$$
In terms of Lemma 2.1, graph $K_t\vee(K_{n-2t-3}\cup K_3\cup tK_1)$ has no perfect $k$-matching. Lemma 2.7 is verified. \hfill $\Box$

\medskip

Let $M$ denote a real symmetric matrix whose rows and columns are indexed by the set $\mathcal{N}=\{1,2,\ldots,n\}$. For a given partition $\pi:\mathcal{N}=\mathcal{N}_1\cup\mathcal{N}_2\cup\cdots\cup\mathcal{N}_t$
of the index set $\mathcal{N}$, the matrix $M$ can be denoted in block form by
\begin{align*}
M=\left(
  \begin{array}{cccc}
    M_{11} & M_{12} & \cdots & M_{1t}\\
    M_{21} & M_{22} & \cdots & M_{2t}\\
    \vdots & \vdots & \ddots & \vdots\\
    M_{t1} & M_{t2} & \cdots & M_{tt}\\
  \end{array}
\right),
\end{align*}
where $M_{ij}$ denotes the $n_i\times n_j$ matrix for $1\leq i,j\leq t$. We denote by $m_{ij}$ the average row sum of $M_{ij}$ for $1\leq i,j\leq t$. Then the matrix $M_{\pi}=(m_{ij})_{t\times t}$ is said to be the
quotient matrix of $M$ with respect to the partition $\pi$. The partition $\pi$ is called equitable if every $M_{ij}$ has constant row sum for $1\leq i,j\leq t$.

\medskip

\noindent{\textbf{Lemma 2.8}} (Brouwer and Haemers \cite{BH}, You, Yang, So and Xi \cite{YYSX}). Let $M$ be a real $n\times n$ matrix with an equitable partition $\pi$, and let $M_{\pi}$ be the corresponding
quotient matrix. Then the eigenvalues of $M_{\pi}$ are eigenvalues of $M$. Furthermore, if $M$ is nonnegative and irreducible, then the largest eigenvalues of $M$ and $M_{\pi}$ are equal.

\section{The proof of Theorem 1.2}

\noindent{\it Proof of Theorem 1.2.} Suppose, to the contrary, that a $t$-connected graph $G$ has no perfect $k$-matching, where $t$ and $k$ are positive integers with $k\equiv1$ (mod 2). By Lemma 2.1,
there exists a subset $S$ of $V(G)$ satisfying $odd(G-S)+k\cdot i(G-S)\geq k|S|+1$. Notice that $n$ is even and $k$ is odd. Then $odd(G-S)+k\cdot i(G-S)$ and $k|S|$ possess the same parity, and so
$odd(G-S)+k\cdot i(G-S)\geq k|S|+2$. We first prove the following claims.

\medskip

\noindent{\bf Claim 1.} $S\neq\emptyset$.

\noindent{\it Proof.} Assume that $S=\emptyset$. Since $n$ is even and $G$ is $t$-connected, we obtain $0=odd(G)+k\cdot i(G)=odd(G-S)+k\cdot i(G-S)\geq k|S|+2=2$. Thus, we get a contradiction. Hence,
$S\neq\emptyset$. Claim 1 is proved. \hfill $\Box$

\medskip

\noindent{\bf Claim 2.} $|S|\geq t$.

\noindent{\it Proof.} Assume that $|S|\leq t-1$. In terms of Claim 1, $k\geq1$, $n\equiv0$ (mod 2) and $G$ being $t$-connected, we deduce
\begin{align*}
1\geq&odd(G-S)+i(G-S)\\
\geq&\frac{1}{k}\cdot odd(G-S)+i(G-S)\\
=&\frac{1}{k}(odd(G-S)+k\cdot i(G-S))\\
\geq&\frac{1}{k}(k|S|+2)\\
\geq&\frac{1}{k}(k+2)\\
>&1,
\end{align*}
which is a contradiction. Therefore, $|S|\geq t$. Claim 2 is proved. \hfill $\Box$

\medskip

\noindent{\bf Claim 3.} $n\geq2|S|+2$.

\noindent{\it Proof.} For $i(G-S)\geq|S|+1$, we have $n\geq|S|+3\cdot odd(G-S)+i(G-S)\geq2|S|+1$. Combining this with $n\equiv0$ (mod 2), we infer $n\geq2|S|+2$.

For $i(G-S)\leq|S|$, we obtain
\begin{align*}
n\geq&|S|+3\cdot odd(G-S)+i(G-S)\\
\geq&|S|+odd(G-S)+i(G-S)\\
=&|S|+(odd(G-S)+k\cdot i(G-S))-(k-1)i(G-S)\\
\geq&|S|+(k|S|+2)-(k-1)|S|\\
=&2|S|+2.
\end{align*}
This completes the proof of Claim 3. \hfill $\Box$

\medskip

Let $|S|=s$, $i(G-S)=i$ and $odd(G-S)=q$. Then we possess $q+ki\geq ks+2$. Let $G_1=K_s\vee(K_{n-2s-1}\cup(s+1)K_1)$.

\medskip

\noindent{\bf Claim 4.} $\rho(G)\leq\rho(G_1)$ with equality if and only if $G=G_1$.

\noindent{\it Proof.} For $i\geq s+1$, it is obvious that $G$ is a spanning subgraph of $G_1=K_s\vee(K_{n-2s-1}\cup(s+1)K_1)$. According to Lemma 2.3, we deduce
$$
\rho(G)\leq\rho(G_1),
$$
with equality if and only if $G=G_1$.

For $i\leq s$, we possess $q\geq k(s-i)+2\geq2$. Obviously, $G$ is a spanning subgraph of $K_s\vee(K_{n_1}\cup K_{n_2}\cup\cdots\cup K_{n_q}\cup iK_1)$, where $n_1\geq n_2\geq\cdots\geq n_q\geq3>1$ are
odd integers with $s+i+\sum\limits_{i=1}^{q}n_i=n$. Using Lemma 2.3, we get
\begin{align}\label{eq:3.1}
\rho(G)\leq\rho(K_s\vee(K_{n_1}\cup K_{n_2}\cup\cdots\cup K_{n_q}\cup iK_1)),
\end{align}
with equality following if and only if $G=K_s\vee(K_{n_1}\cup K_{n_2}\cup\cdots\cup K_{n_q}\cup iK_1)$. Recall that $n_1\geq n_2\geq\cdots\geq n_q\geq3>1$. By virtue of Lemma 2.4, we conclude
\begin{align}\label{eq:3.2}
\rho(K_s\vee(K_{n_1}\cup K_{n_2}\cup\cdots\cup K_{n_q}\cup iK_1))<\rho(K_s\vee(K_{n-s-i-q+1}\cup(i+q-1)K_1)).
\end{align}
Recall that $k\geq1$, $q\geq2$ and $q+ki\geq ks+2$. Then $1+k(i+q-1)=q+ki+(k-1)(q-1)\geq ks+2$, which yields $i+q-1\geq s+1$. Obviously, $K_s\vee(K_{n-s-i-q+1}\cup(i+q-1)K_1)$ is a spanning subgraph of
$G_1=K_s\vee(K_{n-2s-1}\cup(s+1)K_1)$. From Lemma 2.3, we obtain
\begin{align}\label{eq:3.3}
\rho(K_s\vee(K_{n-s-i-q+1}\cup(i+q-1)K_1))\leq\rho(G_1),
\end{align}
where the equality holds if and only if $i+q-1=s+1$. It follows from \eqref{eq:3.1}, \eqref{eq:3.2} and \eqref{eq:3.3} that
$$
\rho(G)<\rho(G_1).
$$
This completes the proof of Claim 4. \hfill $\Box$

Recall that $G_1=K_s\vee(K_{n-2s-1}\cup(s+1)K_1)$. If $s=t$, then we possess $G_1=K_t\vee(K_{n-2t-1}\cup(t+1)K_1)$ and $\rho(G_1)=\rho(K_t\vee(K_{n-2t-1}\cup(t+1)K_1))$. Together with Claim 4, we get
$$
\rho(G)\leq\rho(K_t\vee(K_{n-2t-1}\cup(t+1)K_1)),
$$
where the equality follows if and only if $G=K_t\vee(K_{n-2t-1}\cup(t+1)K_1)$. In terms of Lemma 2.5, $K_t\vee(K_{n-2t-1}\cup(t+1)K_1)$ has no perfect $k$-matching. Thus, we obtain a contradiction. Next,
we shall consider $s\geq t+1$.

For the partition $V(G_1)=V(K_s)\cup V((s+1)K_1)\cup V(K_{n-2s-1})$, the quotient matrix of $A(G_1)$ is
\begin{align*}
B_1=\left(
  \begin{array}{ccc}
  s-1 & s+1 & n-2s-1\\
  s & 0 & 0\\
  s & 0 & n-2s-2\\
  \end{array}
\right).
\end{align*}
The characteristic polynomial of $B_1$ is given by
\begin{align*}
f_{B_1}(x)=x^{3}+(s+3-n)x^{2}+(2-n-s^{2})x-2s^{3}+(n-4)s^{2}+(n-2)s.
\end{align*}
According to Lemma 2.8 and the equitable partition $V(G_1)=V(K_s)\cup V((s+1)K_1)\cup V(K_{n-2s-1})$, $\rho(G_1)$ is the largest root of $f_{B_1}(x)=0$. Namely, $f_{B_1}(\rho(G_1))=0$.

Let $G_*=K_t\vee(K_{n-2t-1}\cup(t+1)K_1)$. The quotient matrix of $A(G_*)$ in view of the partition $V(G_*)=V(K_t)\cup V((t+1)K_1)\cup V(K_{n-2t-1})$ is written as
\begin{align*}
B_*=\left(
  \begin{array}{ccc}
  t-1 & t+1 & n-2t-1\\
  t & 0 & 0\\
  t & 0 & n-2t-2\\
  \end{array}
\right).
\end{align*}
The characteristic polynomial of $B_*$ equals
\begin{align*}
f_{B_*}(x)=&x^{3}+(t+3-n)x^{2}+(2-n-t^{2})x-2t^{3}+(n-4)t^{2}+(n-2)t.
\end{align*}
Based on Lemma 2.8 and the equitable partition $V(G_*)=V(K_t)\cup V((t+1)K_1)\cup V(K_{n-2t-1})$, $\rho(G_*)$ is the largest root of $f_{B_*}(x)=0$. That is, $f_{B_*}(\rho(G_*))=0$.

Let $\rho=\rho(G_*)$. Notice that $G_*=K_t\vee(K_{n-2t-1}\cup(t+1)K_1)$ contains $K_{n-t-1}$ as its proper subgraph. Applying Lemma 2.3, we conclude
\begin{align}\label{eq:3.4}
\rho=\rho(G_*)>\rho(K_{n-t-1})=n-t-2.
\end{align}

By a direct calculation, we obtain
\begin{align}\label{eq:3.5}
f_{B_1}(x)-f_{B_*}(x)=(x-t)g(x),
\end{align}
where $g(x)=x^{2}-(s+t)x-2s^{2}+(n-2t-4)s-2t^{2}+(n-4)t+n-2$. Recall that $\rho=\rho(G_*)>n-t-2$ (see \eqref{eq:3.4}). In what follows, we shall prove that $g(x)>0$ for $x\geq n-t-2$.

The symmetry axis of $g(x)$ is $x=\frac{s+t}{2}<s<n-t-2$, where the last two inequalities occur from the fact $s\geq t+1$ and $n\geq2s+2$ (see Claim 3). This implies that $g(x)$ is increasing in the
interval $[n-t-2,+\infty)$. For $x\geq n-t-2$, we obtain
\begin{align*}
g(x)\geq&g(n-t-2)\\
=&(n-t-2)^{2}-(s+t)(n-t-2)-2s^{2}+(n-2t-4)s-2t^{2}+(n-4)t+n-2\\
=&-2s^{2}-(t+2)s+n^{2}-(2t+3)n+2t+2\\
\geq&-2\Big(\frac{n-2}{2}\Big)^{2}-(t+2)\Big(\frac{n-2}{2}\Big)+n^{2}-(2t+3)n+2t+2 \ \ \ \ \ (\mbox{since} \ s\leq\frac{n-2}{2})\\
=&\frac{1}{2}(n^{2}-(5t+4)n+6t+4)\\
\geq&\frac{1}{2}((5t+3)^{2}-(5t+4)(5t+3)+6t+4) \ \ \ \ \ (\mbox{since} \ n\geq5t+3)\\
=&\frac{t+1}{2}\\
>&0.
\end{align*}
Combining this with \eqref{eq:3.5} and $x\geq t+1$, we infer $f_{B_1}(x)>f_{B_*}(x)$ for $x\geq n-t-2$, which leads to $\rho(G_1)<\rho(G_*)$. Together with Claim 4, we conclude
$$
\rho(G)\leq\rho(G_1)<\rho(G_*)=\rho(K_t\vee(K_{n-2t-1}\cup(t+1)K_1)),
$$
which contradicts $\rho(G)\geq\rho(K_t\vee(K_{n-2t-1}\cup(t+1)K_1))$. This completes the proof of Theorem 1.2. \hfill $\Box$

\section{The proof of Theorem 1.4}

\noindent{\it Proof of Theorem 1.4.} Suppose, to the contrary, that a $t$-connected graph $G$ with a fractional perfect matching has no perfect $k$-matching, where $t$ and $k$ are positive integers with
$k\equiv1$ (mod 2). Based on Lemma 2.1, there exists a subset $S$ of $V(G)$ such that $odd(G-S)+k\cdot i(G-S)\geq k|S|+1$. Since $n$ is even and $k$ is odd, $odd(G-S)+k\cdot i(G-S)$ and $k|S|$ have the
same parity. Hence, we deduce $odd(G-S)+k\cdot i(G-S)\geq k|S|+2$. Obviously, $|S|\geq t$ due to Claim 2.

Let $|S|=s$, $i(G-S)=i$ and $odd(G-S)=q$. Then $q+ki\geq ks+2$. Since $G$ has a fractional perfect matching, it follows from Lemma 2.2 that $i=i(G-S)\leq s$ for every subset $S$ of $V(G)$ with $s\geq t$.
Thus, we obtain $q\geq k(s-i)+2\geq2$. Obviously, $G$ is a spanning subgraph of $G_1=K_s\vee(K_{n_1}\cup K_{n_2}\cup\cdots\cup K_{n_q}\cup iK_1)$, where $n_1\geq n_2\geq\cdots\geq n_q\geq3>1$ are odd
integers with $s+i+\sum\limits_{i=1}^{q}n_i=n$. By Lemma 2.3, we possess
\begin{align}\label{eq:4.1}
\rho(G)\leq\rho(G_1),
\end{align}
with equality if and only if $G=G_1$.

Let $G_2=K_s\vee(K_{n-s-i-q-1}\cup K_3\cup(i+q-2)K_1)$. According to Lemma 2.6, we get
\begin{align}\label{eq:4.2}
\rho(G_1)\leq\rho(G_2),
\end{align}
with equality if and only if $q=2$ and $(n_1,n_2)=(n-s-i-3,3)$. Recall that $k\geq1$, $q\geq2$ and $q+ki\geq ks+2$. Then $2+k(i+q-2)=q+ki+(k-1)(q-2)\geq ks+2$, and so $i+q-2\geq s$.

Let $G_3=K_s\vee(K_{n-2s-3}\cup K_3\cup sK_1)$, where $n\geq2s+6$. It is obvious that $G_2$ is a spanning subgraph of $G_3$. Based on Lemma 2.3, we obtain
\begin{align}\label{eq:4.3}
\rho(G_2)\leq\rho(G_3),
\end{align}
where the equality follows if and only if $G_2=G_3$.

If $s=t$, then $G_3=K_t\vee(K_{n-2t-3}\cup K_3\cup tK_1)$. Together with \eqref{eq:4.1}, \eqref{eq:4.2} and \eqref{eq:4.3}, we deduce
$$
\rho(G)\leq\rho(K_t\vee(K_{n-2t-3}\cup K_3\cup tK_1)),
$$
where the equality occurs if and only if $G=K_t\vee(K_{n-2t-3}\cup K_3\cup tK_1)$. In view of Lemma 2.7, graph $K_t\vee(K_{n-2t-3}\cup K_3\cup tK_1)$ has no perfect $k$-matching. Thus, we obtain a
contradiction. In what follows, we shall consider $s\geq t+1$.

Recall that $G_3=K_s\vee(K_{n-2s-3}\cup K_3\cup sK_1)$. The quotient matrix of $A(G_3)$ based on the partition $V(G_3)=V(K_s)\cup V(K_3)\cup V(sK_1)\cup V(K_{n-2s-3})$ is given by
\begin{align*}
B_3=\left(
  \begin{array}{cccc}
  s-1 & 3 & s & n-2s-3\\
  s & 2 & 0 & 0\\
  s & 0 & 0 & 0\\
  s & 0 & 0 & n-2s-4\\
  \end{array}
\right).
\end{align*}
Based on a computation, the characteristic polynomial of $B_3$ equals
\begin{align*}
\varphi_{B_3}(x)=&x^{4}+(s+3-n)x^{3}+(n-s^{2}-4s-6)x^{2}\\
&+(s^{2}n+3sn+2n-2s^{3}-8s^{2}-14s-8)x-2s^{2}n+4s^{3}+8s^{2}.
\end{align*}
In view of Lemma 2.8 and the equitable partition $V(G_3)=V(K_s)\cup V(K_3)\cup V(sK_1)\cup V(K_{n-2s-3})$, $\rho(G_3)$ is the largest root of $\varphi_{B_3}(x)=0$. That is to say,
$\varphi_{B_3}(\rho(G_3))=0$.

Let $G_*=K_t\vee(K_{n-2t-3}\cup K_3\cup tK_1)$. The quotient matrix of $A(G_*)$ corresponding to the partition $V(G_*)=V(K_t)\cup V(K_3)\cup V(tK_1)\cup V(K_{n-2t-3})$ is
\begin{align*}
B_*=\left(
  \begin{array}{cccc}
  t-1 & 3 & t & n-2t-3\\
  t & 2 & 0 & 0\\
  t & 0 & 0 & 0\\
  t & 0 & 0 & n-2t-4\\
  \end{array}
\right),
\end{align*}
and its characteristic polynomial is
\begin{align*}
\varphi_{B_*}(x)=&x^{4}+(t+3-n)x^{3}+(n-t^{2}-4t-6)x^{2}\\
&+(t^{2}n+3tn+2n-2t^{3}-8t^{2}-14t-8)x-2t^{2}n+4t^{3}+8t^{2}.
\end{align*}
By virtue of Lemma 2.8 and the equitable partition $V(G_*)=V(K_t)\cup V(K_3)\cup V(tK_1)\cup V(K_{n-2t-3})$, the largest root of $\varphi_{B_*}(x)=0$ is equal to $\rho(G_*)$. Namely,
$\varphi_{B_*}(\rho(G_*))=0$.

By a computation, we obtain
\begin{align}\label{eq:4.4}
\varphi_{B_3}(x)-\varphi_{B_*}(x)=(s-t)\psi(x),
\end{align}
where $\psi(x)=x^{3}-(s+t+4)x^{2}+(sn+tn+3n-2s^{2}-(2t+8)s-2t^{2}-8t-14)x+4s^{2}+(4t+8-2n)s-2tn+4t^{2}+8t$. Since $G_*$ contains $K_{n-t-3}$ as its proper subgraph, it follows from Lemma 2.3 that
$\rho(G_*)>\rho(K_{n-t-3})=n-t-4$. In what follows, we shall verify that $\psi(x)>0$ for $x\geq n-t-4$.

By a direct calculation, we have $\psi'(x)=3x^{2}-2(s+t+4)x+sn+tn+3n-2s^{2}-(2t+8)s-2t^{2}-8t-14$ and the symmetry axis of $\psi'(x)$ is $x=\frac{s+t+4}{3}$. According to $s\geq t+1$ and $n\geq2s+6$,
we deduce
$$
\frac{s+t+4}{3}<n-t-4,
$$
which yields that $\psi'(x)$ is increasing in the interval $[n-t-4,+\infty)$. When $x\geq n-t-4$, it follows from $s\leq\frac{n-6}{2}$ and $n\geq5t+7$ that
\begin{align*}
\psi'(x)\geq&\psi'(n-t-4)\\
=&3(n-t-4)^{2}-2(s+t+4)(n-t-4)+sn+tn+3n\\
&-2s^{2}-(2t+8)s-2t^{2}-8t-14\\
=&-2s^{2}-ns+3n^{2}-(7t+29)n+3t^{2}+32t+66\\
\geq&-2\Big(\frac{n-6}{2}\Big)^{2}-n\Big(\frac{n-6}{2}\Big)+3n^{2}-(7t+29)n+3t^{2}+32t+66\\
=&2n^{2}-(7t+20)n+3t^{2}+32t+48\\
\geq&2(5t+7)^{2}-(7t+20)(5t+7)+3t^{2}+32t+48\\
=&18t^{2}+23t+6\\
>&0,
\end{align*}
which implies that $\psi(x)$ is increasing in the interval $[n-t-4,+\infty)$. For $x\geq n-t-4$, we obtain
\begin{align}\label{eq:4.5}
\psi(x)\geq&\psi(n-t-4)\nonumber\\
=&(n-t-4)^{3}-(s+t+4)(n-t-4)^{2}\nonumber\\
&+(sn+tn+3n-2s^{2}-(2t+8)s-2t^{2}-8t-14)(n-t-4)\nonumber\\
&+4s^{2}+(4t+8-2n)s-2tn+4t^{2}+8t\nonumber\\
=&-(2n-2t-12)s^{2}-((t+6)n-t^{2}-12t-24)s\nonumber\\
&+n^{3}-(3t+13)n^{2}+(2t^{2}+23t+54)n-4t^{2}-42t-72\nonumber\\
\geq&-(2n-2t-12)\Big(\frac{n-6}{2}\Big)^{2}-((t+6)n-t^{2}-12t-24)\Big(\frac{n-6}{2}\Big)\nonumber\\
&+n^{3}-(3t+13)n^{2}+(2t^{2}+23t+54)n-4t^{2}-42t-72\nonumber\\
& \ \ \ \ \ \ \ \ \ \ \Big(\mbox{since} \ t+1\leq s\leq\frac{n-6}{2} \ \mbox{and} \ n\geq5t+7\Big)\nonumber\\
=&\frac{1}{2}(n^{3}-(6t+14)n^{2}+(5t^{2}+52t+60)n-14t^{2}-120t-72)\nonumber\\
=&\frac{1}{2}T(n),
\end{align}
where $T(n)=n^{3}-(6t+14)n^{2}+(5t^{2}+52t+60)n-14t^{2}-120t-72$. By a direct calculation, we get $T'(n)=3n^{2}-2(6t+14)n+5t^{2}+52t+60$. Notice that
$$
\frac{6t+14}{3}<5t+7\leq n.
$$
Thus, we possess
\begin{align*}
T'(n)\geq&T'(5t+7)\\
=&3(5t+7)^{2}-2(6t+14)(5t+7)+5t^{2}+52t+60\\
=&20t^{2}+38t+11\\
>&0.
\end{align*}
This implies that $T(n)$ is increasing in the interval $[5t+7,+\infty)$, and so
\begin{align*}
T(n)\geq&T(5t+7)\\
=&(5t+7)^{3}-(6t+14)(5t+7)^{2}\\
&+(5t^{2}+52t+60)(5t+7)-14t^{2}-120t-72\\
=&36t^{2}+5t+5\\
>&0.
\end{align*}
Combining this with \eqref{eq:4.5}, we obtain that $\psi(x)>0$ for $x\geq n-t-4$. Together with \eqref{eq:4.4} and $s\geq t+1$, we possess $\varphi_{B_3}(x)>\varphi_{B_*}(x)$ for $x\geq n-t-4$,
which leads to $\rho(G_3)<\rho(G_*)$. Combining this with \eqref{eq:4.1}, \eqref{eq:4.2} and \eqref{eq:4.3}, we conclude
$$
\rho(G)\leq\rho(G_1)\leq\rho(G_2)\leq\rho(G_3)<\rho(G_*)=\rho(K_t\vee(K_{n-2t-3}\cup K_3\cup tK_1)),
$$
which contradicts $\rho(G)\geq\rho(K_t\vee(K_{n-2t-3}\cup K_3\cup tK_1))$. This completes the proof of Theorem 1.4. \hfill $\Box$

\section*{Declaration of competing interest}

The authors declare that they have no known competing financial interests or personal relationships that could have appeared to influence the work reported in this paper.

\section*{Data availability}

No data was used for the research described in the article.

\section*{Acknowledgments}

This work was supported by the Natural Science Foundation of Jiangsu Province (Grant No. BK20241949). Project ZR2023MA078 supported by Shandong Provincial Natural Science Foundation.

\end{document}